\documentclass[12pt]{amsart}
\usepackage{amssymb, epsf, colordvi}
\usepackage{multirow}
\usepackage{latexsym}
\usepackage{graphics}

\newif\ifEPSF
\EPSFtrue %
\numberwithin{equation}{section}

\newtheorem{thm}{Theorem}
\numberwithin{thm}{section}

\newtheorem{example}[thm]{Example}
\newtheorem{remark}[thm]{Remark}
\newtheorem{definition}[thm]{Definition}

\newcounter{FNC}[page]
\def\newfootnote#1{{\addtocounter{FNC}{2}$^\fnsymbol{FNC}$%
     \let\thefootnote\relax\footnotetext{$^\fnsymbol{FNC}$#1}}}

\copyrightinfo{}{}

\newcommand\trop{\operatorname{trop}}
\title{Real Zeuthen numbers for two lines}

\author{Benoit Bertrand}
\address{Section de Math\'ematiques\\
        Universit\'e de Gen\`eve\\
        Case postale 64\\
        1211 Gen\`eve 4\\
        Suisse}
\email{benoit.bertrand@math.unige.ch}
\urladdr{http://www.unige.ch/math/folks/bertrand}

\thanks{The author  was partially supported by the European research
  network IHP-RAAG contract HPRN-CT-2001-00271 and part of the work
  was done at Max Planck Institute for Mathematics in Bonn}

\dedicatory{Dedicated to the memory of Felice Ronga.}

\newtheorem{bbthm}{Theorem}

\newcommand\CQFD{\hfill $\Box$ \newline}

\newcommand\ZZ{{\mathbb{Z}}}
\newcommand\PP{{\mathbb{P}}}
\newcommand\RR{{\mathbb{R}}}
\newcommand\CC{{\mathbb{C}}}

\newcommand\mult{\operatorname{mult}}

\begin{document}

\pagestyle{plain}
\setcounter{page}{1}

\begin{abstract}
  Given three natural numbers $k,l,d$ such that $k+l=d(d+3)/2$, the
  Zeuthen number $N_{d}(l)$ is the number of nonsingular complex
  algebraic curves of degree $d$ passing through $k$ points and
  tangent to $l$ lines in $\PP^2$.  It does not depend on the generic
  configuration $C$ of points and lines chosen.  If the points and
  lines are real, the corresponding number $N_{d}^\RR(l,C)$ of
  real curves usually depends on the configuration chosen.  We use
  Mikhalkin's tropical correspondence theorem to prove that for two
  lines the real Zeuthen problem is
  maximal: there exists a configuration $C$ such that
  $N_{d}^\RR(2,C)=N_{d}(2)$. The correspondence theorem
  reduces the computation to counting certain lattice paths with
  multiplicities.
\end{abstract}
\maketitle

\tableofcontents



\section*{Introduction}

Given $l$ lines and $k = d(d+3)/2 - l$ points in $\CC P^2$, how many
nonsingular complex algebraic curves of degree $d$ pass through the
$k$ points and are tangent to the $l$ lines? This is a particular
instance of the Zeuthen problem.  For generic configurations of points
and lines there are finitely many solutions to the problem and we call
Zeuthen number the number $N_{d}(l)$ of solutions. Here we consider
the corresponding question for real data: assume that the points and
the lines are real, how many degree $d$ real curves pass through the
$k$ points and are tangent to the $l$ lines? In other words, what
values can take the real Zeuthen number ${N_d^\RR} (l,C)$ of real
solutions?  This number usually depends on the configuration $C$, and,
clearly, the (invariant) number of complex solutions is an upper
bound.
  Whether there exists a generic configuration for which all the
solutions are real is a natural and classical question in real
enumerative geometry.  It is said that the problem is maximal if such
a configuration exists.  For $l=1$ it was shown by F.  Ronga
\cite{Ron00} that the Zeuthen problem is maximal in the above sense
(i.e. all the curves can be real).  In this
article we show that the problem for $2$ lines is also maximal.
\begin{bbthm}\label{felice}
For any integer $d \ge 2$ there exists a configuration $C$ of $2$ real
lines and $d(d+3)/2 -2$ real points such that all the degree $d$ curves passing
through the points and tangent to the lines are real:
$$N^\RR_{d} (2,C) = N_{d} (2)$$
\end{bbthm}
The techniques we use are those developed by Mikhalkin in
\cite{mikh05}.  The statement is proved using correspondence theorems
to tropicalize the problem and the lattice path algorithm to count the number
of tropical curves. I would like to thank G. Mikhalkin and E. Shustin
for useful discussions.



\section{Tropical tangency and algorithm}


\subsection{Tropical tangency}

\begin{figure} 
  \resizebox {10 cm}{!}{\includegraphics{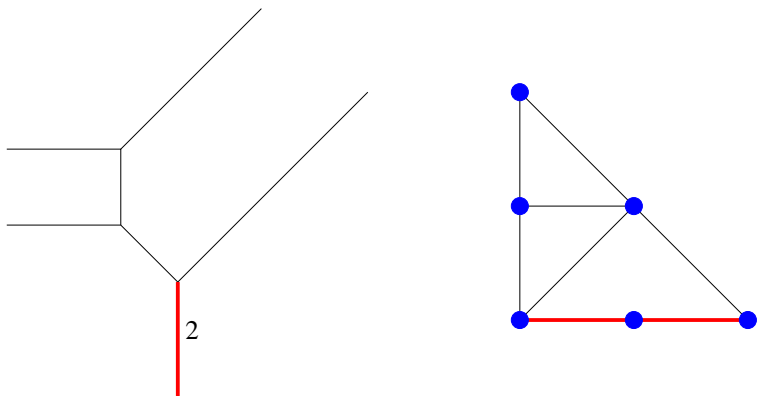}} 
\caption{Tropical conic tangent to the horizontal axis and dual
  triangulation of its Newton polygon.}
\label{tang}
\end{figure}

We use here the terminology introduced in \cite{mikh05} without
recalling definitions and basic properties of tropical curves.
General tropical tangency has been discussed in Mikhalkin's note
\cite{MikhXN}. Here we consider a simple special case of tangency with
respect to toric divisors i.e. in our projective plane setting, the
three axes which are the closures of one dimensional orbits of the
torus action (See also \cite{GatMar07}). Recall that a plane tropical
curve is dual to a subdivision of its Newton polygon such that the
weight of an edge of the curve is equal to lattice length its dual
edge. A simple tangency to a coordinate axis in the algebraic world
translates to a weight $2$ unbounded edge in tropical world.
\begin{definition}
  We say that a tropical curve $C$ with Newton polygon $\Delta$ is
  tangent to an axis corresponding to the edge $\delta \subset \Delta$
  if in the dual subdivision there is a $1$-simplex of length greater
  or equal to $2$ included in $\delta$ (See Figure~\ref{tang}). 
\end{definition}

We only consider here simple tangencies (i.e. the weight of the
corresponding edge is 2).

Under the degeneration process described in \cite{mikh05} Section~6.1
a family of curves tangent to some axes and passing through the right
number of points has a tropical limit curve which is tangent (in the
above sense) to the axes.

All configurations considered here will be generic.  The multiplicity
$\mult(A)$ of a tropical curve $A$ and the multiplicity $\mult_\RR(A)$
of a real tropical curve 
 is defined as in
\cite{mikh05}, Definition~4.15 and~7.10. The proof of
Theorem~\ref{felice} we give below will use the so called lattice path
algorithm of counting tropical curves with multiplicity. The
multiplicity of a path is explained below in subsection~\ref{lpath}.

\begin{bbthm}\label{corres}
Let $l \le 3$ be a positive integer.
The number of complex algebraic curves of degree $d$ passing through
$d(d+3)/2 -l$ points and tangent to $l$ lines is equal to the number
of degree $d$ tropical curves through $d(d+3)/2 -l$ points and tangent
to $l$ chosen axes counted with multiplicities.
\end{bbthm}

\begin{bbthm}\label{real-corres}
Let $l \le 3$ be a positive integer.  Let $T$ be a configuration of
$d(d+3)/2 -l$ 
 points. Let $N^\RR_{\trop,d,T}$ be the number of
degree $d$ real tropical curves counted with multiplicities passing
through $C$ and tangent to $l$ chosen axes. Then there exists a
configuration $C$ of $d(d+3)/2 -l$ real points and $l$ real lines such
that $N^\RR_{d}(l,C)= N^\RR_{\trop,d,T}$.
\end{bbthm}

The first statement can be extracted from \cite{MikhXN} and
\cite{mikh05} and the second one can be deduced from the first one
analyzing, as in \cite{mikh05}, which tropical curves actually produce
real algebraic curves (See \cite{mikh05} Section~8 p.~365 and
Propositions~6.17, 6.18, 8.21 and Theorems~1 and 3).  Alternatively one
can also deduce Theorem~\ref{corres} from E. Shustin's Lemma~5.4 and Theorem~5 in
\cite{Shu06}.




\subsection{Lattice path algorithms}\label{lpath}

As in \cite{mikh05} Sections~7.2 and 7.4, placing the points on an line
of negative irrational slope and sufficiently far away from one
another, one obtains algorithms to count the relevant tropical curves.
We describe below in our special cases the slightly modified
 Mikhalkin's algorithms (See \cite{mikh05} Sections~7.2
and 7.4) that take tangencies into account. We refer to \cite{mikh05} for 
a more developed presentation of lattice path algorithms.



Let $\Delta$ be the triangle with vertices $(0,0),(d,0),(0,d)$ (that
is the Newton polygon of a generic curve of degree $d$).  Let $l \le
3$ be a positive integer and $p_1,\dotsc ,p_l$ be integer points on
$l$ different edges of $\Delta$.  Choose a nearly horizontal line $L$
of $\RR^2$ with a very small negative irrational slope and consider
the orthogonal projection $\lambda: \RR^2 \to L \sim \RR$.  Let $p$
(resp. $q$) be the point of $\Delta$ where $\lambda$ reaches it's
minimum (resp. maximum).  Let $\gamma : [0,n] \to \Delta$ be an
increasing (piecewise linear) lattice path from $p$ to $q$ avoiding the $p_i$'s
(i.e. $\gamma([0,n]\cap\ZZ)\subset\Delta\cap\ZZ^2\setminus
\{p_1,\dotsc,p_l\}$ and $\lambda \circ \gamma$ is increasing).  We
define inductively positive and negative multiplicities $\mu_\pm$ of
$\gamma$.  Let $\alpha_+$ (resp. $\alpha_-$) be the path supported by
the upper edge of $\Delta$ (resp. the vertical and lower edges of
$\Delta$).  The multiplicities $\mu_\pm(\alpha_\pm)$ of
$\alpha_\pm$ are $1$.  If the path $\gamma$ is neither $\alpha_+$ nor
$\alpha_-$ then it divides $\Delta$ into two closed regions $\Delta_+$
and $\Delta_-$, where $\Delta_\pm$ contains $\alpha_\pm$. Let $k$ be
the smallest integer such that $\Delta_\pm$ is locally strictly convex at
$\gamma(k)$.
  Consider the path $\gamma^\prime : [0,n-1] \to \Delta$ defined by 
$\gamma^\prime (j)= \gamma(j)$ if $j<k$ and  $\gamma^\prime (j)=
\gamma(j+1)$ if $j\ge k$. We set:

\[
\mu_\pm(\gamma) :=2 area(T) \mu_\pm(\gamma^\prime), 
\]
 
where $T$ is the triangle with vertices
$\gamma(k-1),\gamma(k),\gamma(k+1) $.  The multiplicity $\mu(\gamma)$
is defined by $\mu(\gamma)=\mu_+(\gamma)\cdot \mu_-(\gamma)$. Here is
 the statement in our particular case.

Let $\delta_1, \dotsc ,\delta_3$ be edges of $\Delta$ and $\eta_i$ be
the set of integer points in the relative interior of $\delta_i$.

\begin{bbthm}
Let $l \le 3$ be a positive integer.
The number $N_{d} (l)$ is equal to the sum of the multiplicities the
$\lambda$-increasing lattice paths $[0,d(d+3)/2 - l] \to \Delta$
avoiding $p_1, \dotsc, p_l$ over all $\{p_1, \dotsc, p_l\} \in
\eta_1 \times \cdots \times \eta_l$.
\end{bbthm}


There is a real version of the algorithm to count real
curves. We briefly explain below how to obtain it from the previous one
in our special case and refer to \cite{mikh05} for details.

One now needs to consider 'signs' on intervals $[j-1, j]$ for $j= 1,
\dotsc, 
n $
 and 'phases' of edges which are the images of
those intervals.  Let $s \in \ZZ_2^2$ be a choice of `sign' on each
interval. Suppose $\gamma(j)-\gamma(j-1) = (y_j,x_j)\in\ZZ^2, j=1,
\dotsc, 
 n$.
 Let $S_j$ be the quotient of $\ZZ_2^2$ by the equivalence relation
$(X,Y)\sim (X+x_i \mod 2,Y+y_i \mod 2), (X,Y) \in \ZZ_2^2$.  It
induces a phase $\sigma_j (s) \in S_j$ on each corresponding edge
$[\gamma(j-1), \gamma(j)]$ of the image of $\gamma$ in $\Delta$. If $S
= (S_1, \dotsc, S_n) \in (\ZZ_2^2)^{n}$ is a choice of 'signs' we
denote by $\sigma(S)$ the $n$-tuple of induced phases $\sigma_1 (S_1),
\dotsc, \sigma_n (S_n)$

Let $\sigma = (\sigma_j)_{j=1, \dotsc, n},\sigma_j\in S_j$ be any
choice of phases on the edges of the image of $\gamma$.  The real
multiplicity $\mu^\RR(\gamma,\sigma)$ of a path $\gamma$ equipped with
phases $\sigma$ is given by $\mu^\RR(\gamma,\sigma) =
\mu_+^\RR(\gamma,\sigma)\cdot\mu_-^\RR(\gamma,\sigma)$ where
$\mu_\pm^\RR (\gamma,\sigma)$ are defined inductively as above.
The path $\gamma^\prime$ is obtained as before and we set 
$\mu_\pm^\RR (\alpha_\pm,\sigma)=1$ and


\[
\mu_\pm^\RR(\gamma,\sigma) = a(T) \mu_\pm^\RR(\gamma^\prime,\sigma^\prime), 
\]

where 

\begin{itemize}
\item $a(T)=1$ if the lattice area (twice the usual one) of $T$ is
odd. In this case the phase $\sigma_k^\prime$ is chosen such that the
three classes $\sigma_k$, $\sigma_{k+1}$ and $\sigma_k^\prime$ do not
share a common element.

\item $a(T)=0$ if all sides of $T$ have even lattice length and
  $\sigma_k \neq \sigma_{k+1}$.
\item $a(T)=4$ if  all sides of $T$ have even lattice length and
  $\sigma_k = \sigma_{k+1}$. Then define $\sigma_k^\prime :=
  \sigma_k$.
\item $a(T)=0$ if $\sigma_k$ and $\sigma_{k+1}$ do not have a common
  element.
\item $a(T)=2$ if $T$ has exactly one even side distinct from
  $[\gamma(k-1),\gamma(k+1)]$. Then $\sigma_k^\prime$ is defined by the fact
that it should have a common element with $\sigma_k$ and $\sigma_{k+1}$.
\item if $[\gamma(k-1),\gamma(k+1)]$ is the only even side of $T$ then
  one should consider the path $\gamma^\prime$ with the two choices of
  $\sigma_k^\prime$ satisfying the above condition. (i.e.
  $\mu_\pm^\RR(\gamma,\sigma) = \mu_\pm^\RR(\gamma^\prime,\sigma_1^\prime)
  + \mu_\pm^\RR(\gamma^\prime,\sigma_2^\prime) $).
\end{itemize} 

\begin{bbthm}
Let $l \le 3$ be a positive integer.  For any choice of $S \in
  (\ZZ_2^2)^{d(d+3)/2 - l}$ there exists a configuration of $d(d+3)/2 -
  l$ generic points such that the number of degree $d$ real curves
  among the $N_{d} (l)$ complex ones is equal to the sum of the
  multiplicities $\mu^\RR (\gamma,\sigma(S))$ of the $\lambda$-increasing
  signed lattice paths $(\gamma,\sigma(S)) :[0,d(d+3)/2 - l] \to \Delta$
  avoiding $p_1, \dotsc, p_l$ over all $\{p_1, \dotsc, p_l\} \in
  \delta_1 \times \cdots \times \delta_l$.
\end{bbthm}

\section{Maximality}
\label{Maximality}

\subsection{Case of one line}\label{k1}
The case of one line was studied by Felice Ronga \cite{Ron00} who
proved maximality and completeness of the real Zeuthen problem. We give
below a tropical proof of the maximality for one line.

\begin{figure}
  \begin{tabular}{cc}
    \begin{minipage}{7cm}
      \resizebox {6 cm}{!}{\includegraphics{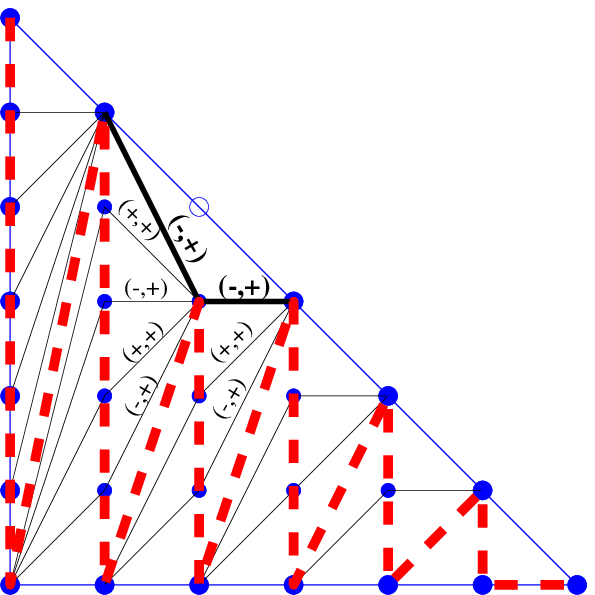}}\\

      Tangency to the line at infinity.

\end{minipage}
&
\begin{minipage}{7cm}

  \resizebox {6.5 cm}{!}{\includegraphics{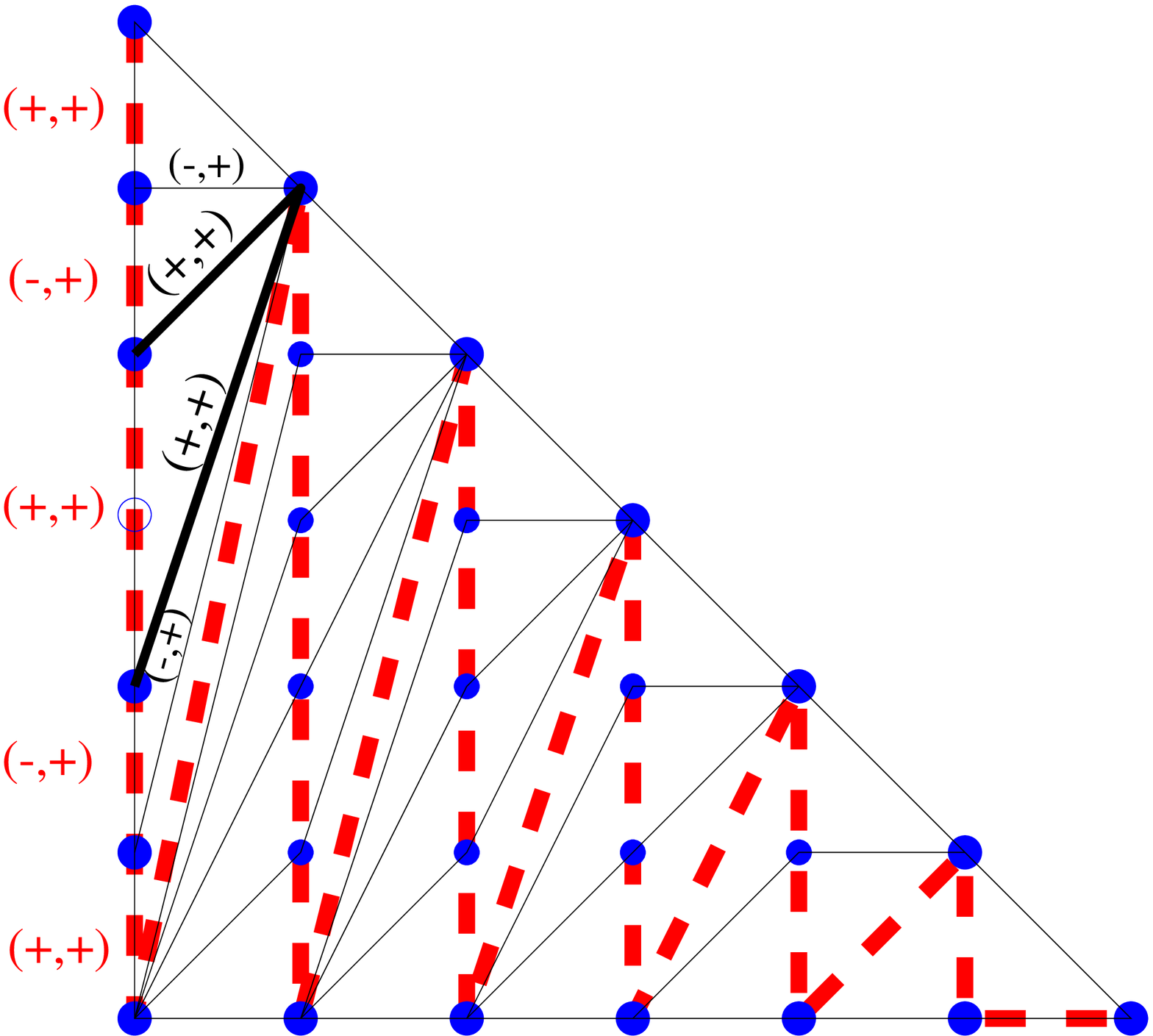}}\\

  Tangency to the vertical axis.
\end{minipage}
\end{tabular}
\caption{Case of one line.}
\label{fig1}
\end{figure}

As we can choose the line to be any of the axes we will consider two
cases: first the line at infinity and then ordinate axis.  To count
the lines tangent to the infinity line we need to count lattice paths
of maximal length avoiding one the interior points $(m,d-m)$ of the
hypothenus $h$ of $\Delta$.  Consider the constant sign sequence
$((+,+), \dotsc, (+,+))$ (we denote by $+$ the zero element of $\ZZ_2$
and by $-$ the nonzero one). We have $d-1$ choices for the point
$(m,d-m)$ and, as shown on Figure~\ref{fig1}, for each one the lattice
path has multiplicity $2$. Thus the total number of curves passing
through our configuration is $2(d-1)$.  When considering the ordinate
axis, the paths must avoid a point on the vertical edge $v$ of
$\Delta$ and we choose the sign sequence $S_{d(d+3)/2 -1}$ which is
constant and equal to $(+,+)$ except for the $d-1$ first terms that
are $(-,+),(+,+),(-,+),(+,+), \dotsc, (-,+),(+,+) $ for $d$ odd and
$(+,+),(-,+),(+,+), \dotsc, (-,+),(+,+) $ for $d$ even. Again as shown
on Figure~\ref{fig1} each path has multiplicity $2$.

\subsection{Main theorem}

We prove that for two lines there is a configuration 
of $d(d+3)/2 - 2$ points and two lines for which all the (nonsingular)
curves of degree $d$ passing through the points and tangent to the
lines are real.  


\proof
We can choose the lines to be the infinity line and the ordinate axis.
Consider the sign sequence $S_{d(d+3)/2 -2}$ of length $d(d+3)/2 - 2$
defined in Section~\ref{k1}. We just need to see that for each pair of
integer points lying respectively in the interior of the edges $h$ and
$v$ we have a lattice path of multiplicity $4$ to prove that
$N^\RR_{d} (2,C)= (2(d-1))^2$.

\begin{figure}[ht]
  \begin{tabular}{cc}
    \begin{minipage}{7cm}
      \resizebox {6.5 cm}{!}{\includegraphics{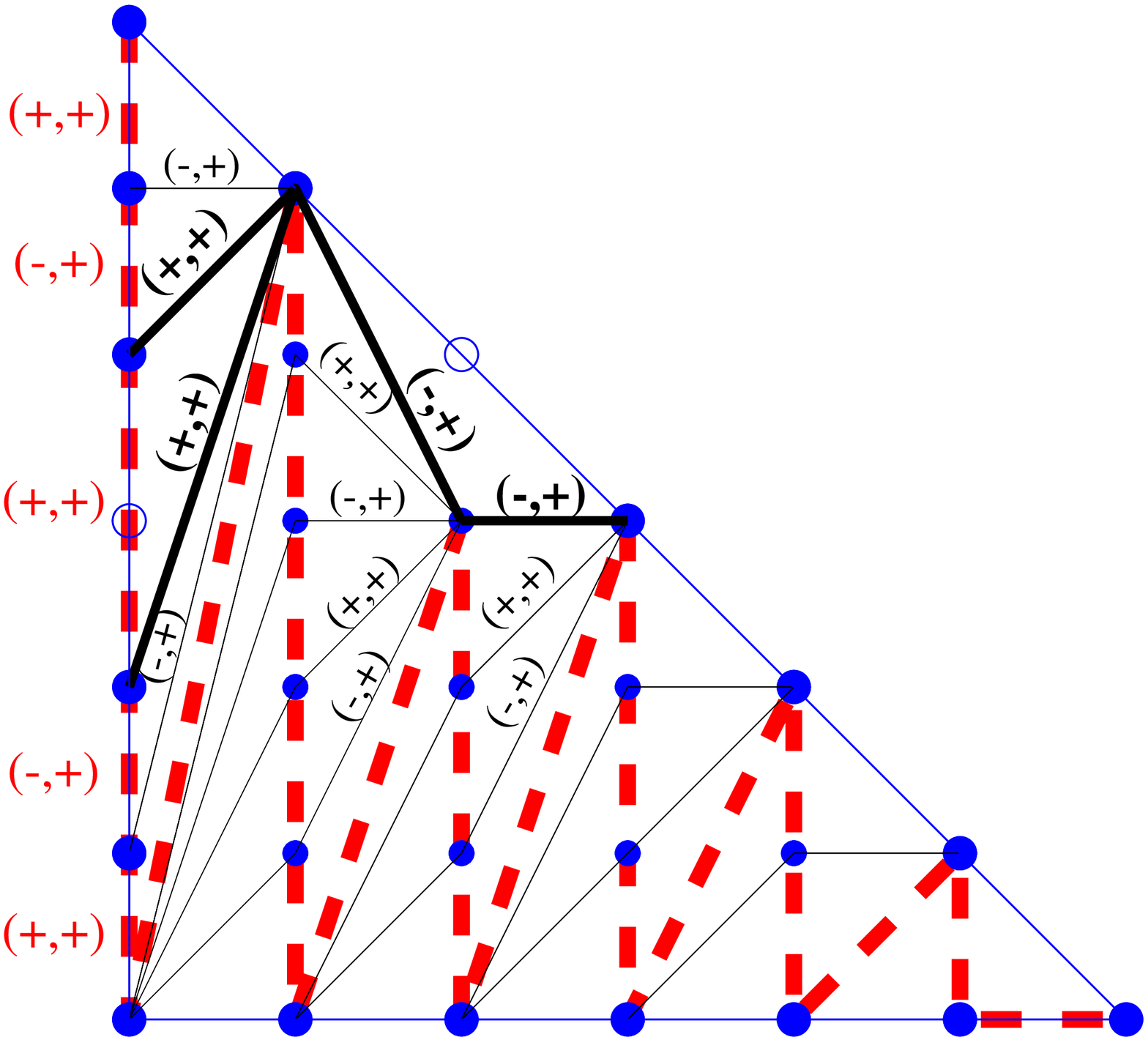}}\\

      Splitting case.

\end{minipage}
&
\begin{minipage}{7cm}

  \resizebox {6 cm}{!}{\includegraphics{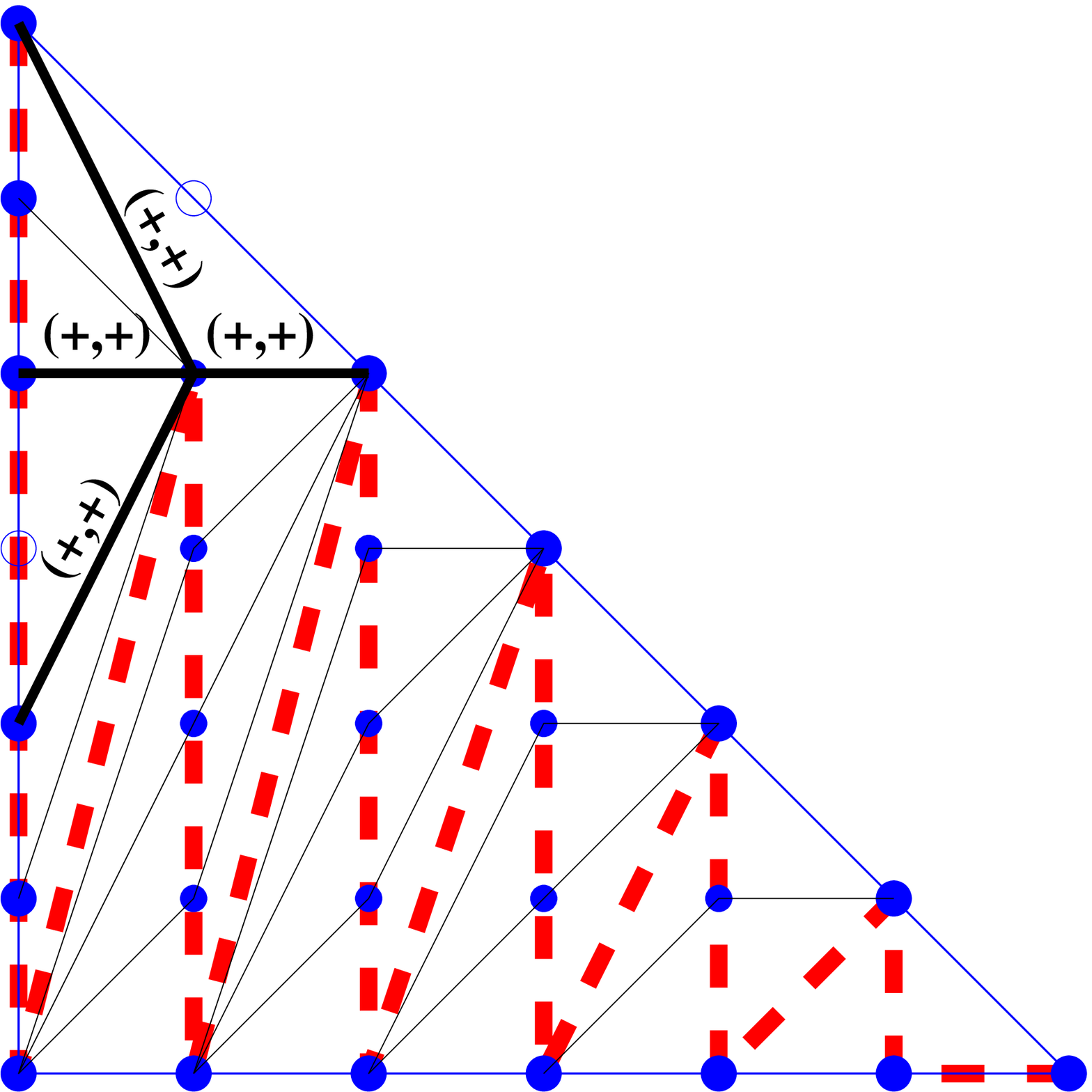}}\\

  Non independent case.
\end{minipage}
\end{tabular}
\caption{The two lines case.}
\label{figure:fig2}
\end{figure}

If the point on $h$ is not $(1,d-1)$ then the picture splits into two
independent parts: the band $[0,1]\times [0,d]$ and the triangle with
vertices $(0,1)$, $(0,d)$ and $(d-1,1)$. Both part where studied for
the case $l=1$ in previous section and contribute each for $2$ in the
multiplicity. (See Figure \ref{figure:fig2})

If the first point is $(1,d-1)$, one just needs to notice that the
signed segments on $v$ behave exactly as if the one of length $2$ was
split in two of length $1$ and the sign distribution was that of
second case of Section~\ref{k1}.  Hence the multiplicity of these
paths is also $4$.(See Figure \ref{figure:fig2}) \CQFD

\def\cprime{$'$} \def\cprime{$'$} \def\cprime{$'$}


\begin{thebibliography}{Ron00}

\bibitem[GM07]{GatMar07}
Andreas Gathmann and Hannah Markwig.
\newblock The {C}aporaso-{H}arris formula and plane relative {G}romov-{W}itten
  invariants in tropical geometry.
\newblock {\em Math. Ann.}, 338(4):845--868, 2007.

\bibitem[Mik]{MikhXN}
Grigory Mikhalkin.
\newblock {Tropical computation of Zeuthen's numbers of toric surfaces}.
\newblock http://www.math.toronto.edu/mikha/XN.ps.

\bibitem[Mik05]{mikh05}
Grigory Mikhalkin.
\newblock Enumerative tropical algebraic geometry in {$\RR^2$}.
\newblock {\em J. Amer. Math. Soc.}, 18(2):313--377 (electronic), 2005.

\bibitem[Ron00]{Ron00}
Felice Ronga.
\newblock A real {R}iemann-{H}urwitz theorem.
\newblock {\em Bol. Soc. Brasil. Mat. (N.S.)}, 31(2):175--187, 2000.

\bibitem[Shu05]{Shu06}
Eugenii Shustin.
\newblock A tropical approach to enumerative geometry.
\newblock {\em Algebra i Analiz}, 17(2):170--214, 2005.

\end{thebibliography}
\end{document}